\documentclass{amsart}
\usepackage[frenchb,english]{babel}
\usepackage{graphicx}
\usepackage{subcaption}
\usepackage{mathabx}

\usepackage{stmaryrd}
\usepackage{appendix}
\usepackage{verbatim}
\usepackage{varwidth}
\usepackage{natbib}

\usepackage[margin=1.2in]{geometry}
\usepackage{multirow}
\newtheorem{theorem}{Theorem}[section]

\theoremstyle{definition}
\newtheorem{definition}[theorem]{Definition}

\theoremstyle{remark}

\numberwithin{equation}{section}
\newcommand{\norm}[1]{\left\lVert#1\right\rVert}
\newcommand{\ie}{\textit{i.e.,\ }}
\newcommand{\eg}{\textit{e.g.,\ }}
\newcommand{\dd}{\textrm{d}}

\newcommand{\DOF}{DoF\,}

\newcommand{\R}{\mathbb{R}}     \newcommand{\Z}{\mathbb{Z}}

   \newcommand{\HH}{\mathcal{H}}
\newcommand{\KK}{\mathcal{K}}

   \newcommand{\SC}{\mathcal{S}}
   \newcommand{\UU}{\mathcal{U}}
   
   \newcommand{\PP}{\mathcal{P}}
   \newcommand{\FF}{\mathcal{F}}

\usepackage{rotating}
\usepackage{breakcites}
\usepackage{lineno}

\begin{document}
\title{
 Symbolic dynamics in a binary asteroid system}
\author{Sara Di Ruzza}
\address{Dipartimento di Matematica Tullio Levi--Civita}
\curraddr{via Trieste, 63, 35121, Padova}
\email{sdiruzza@math.unipd.it}
\thanks{This research is funded by the ERC project 677793 StableChaoticPlanetM}
\author{J\'er\^ome Daquin}
\address{Dipartimento di Matematica Tullio Levi--Civita}
\curraddr{via Trieste, 63, 35121, Padova}
\email{daquin.jerome@gmail.com}
\author{Gabriella Pinzari}
\address{Dipartimento di Matematica Tullio Levi--Civita}
\curraddr{via Trieste, 63, 35121, Padova}
\email{gabriella.pinzari@math.unipd.it}


\date{\today}

\date{}

\dedicatory{}

\begin{abstract}
We highlight the existence of a topological horseshoe arising from a a--priori stable model of the binary asteroid dynamics. The inspection is numerical and uses correctly aligned windows, as described in a recent paper by A. Gierzkiewicz and P. Zgliczy\'nski, combined with a recent analysis of an associated secular problem.
 \end{abstract}
 
\maketitle
\tableofcontents

\section{Purpose of the paper}
\noindent This paper aims to highlight chaos in the secular motions of a binary asteroid system interacting with a planet whose orbit is external to the orbits of the asteroids. These chaotic motions turn to bifurcate from an  a--priori stable configuration, in the sense of~\cite{chierchiaG94}.  We shall not provide rigorous proofs, besides the heuristic arguments that we are going to present in this introduction. In fact, our study will be purely numerical. Moreover, we shall  not implement any algorithm to control machine errors. We are however convinced that our computations are correct thanks to a--posteriori checks that we shall describe in the course of the paper.  \\
 Let us  describe the physical setting. Three point masses constrained on a plane undergo Newtonian attraction. Two of them (the asteroids) have comparable (in fact, equal) mass and, approximately, orbit their common barycentre. The orbit of a much more massive body (the planet) keeps external to the  couple, for a sufficiently long time. We do not assume\footnote{See, \eg~\cite{paezL15} for a study based on a restricted model.} any prescribed trajectory for any of the bodies, but just Newton law as a mutual interaction. We fix a reference frame centred with one of the asteroids and we look at the motions of the other one  and the planet. As no Newtonian interaction can be regarded as dominant -- as, for example, in the cases investigated in~\cite{arnold63, fejoz04, laskarR95, pinzari-th09, chierchiaPi11b} and~\cite{gioLS17, volpiLS18} -- in order to simplify the analysis, we  look at a certain {\it secular} system, obtained, roughly,  averaging out  the proper time of the reference asteroid. This means that  we are assuming that the time scale of the movements of the planet is much longer. Beware that our secular problem  has nothing to do with the  one usually considered in the literature, where the average is performed with respect to {\it two} proper times (\eg~\cite{fejozG2016}). Let us look, for a moment, to the case where the planet  is constrained on a circular trajectory.  In such case, the only observables are the eccentricity and the pericentre of the instantaneous ellipse of the  asteroid. Quantitatively, this system may be described by only two conjugate Hamiltonian coordinates: the angular momentum $G$ (related to the eccentricity) and the pericentre coordinate $g$ of the asteroidal ellipse. There is a limiting situation, which roughly corresponds to  the planet being at infinite distance, where, exploiting results from~\cite{gPi19, gPi20a, gPi20b},  the phase portrait of the system in the plane $(g, G)$ reveals only librational periodic motions. Physically, such motions correspond to the perihelion direction of the asteroidal ellipse affording small oscillations about one equilibrium position, with the ellipse highly eccentric and periodically squeezing to a segment. The movements are accompanied by a change of sense of motion every half--period. The purpose of this paper is to highlight the onset of chaos in the full secular problem, when the planet is far and moves almost circularly.  

\begin{figure}[h!]
\includegraphics[height=9.5cm, width=13cm]{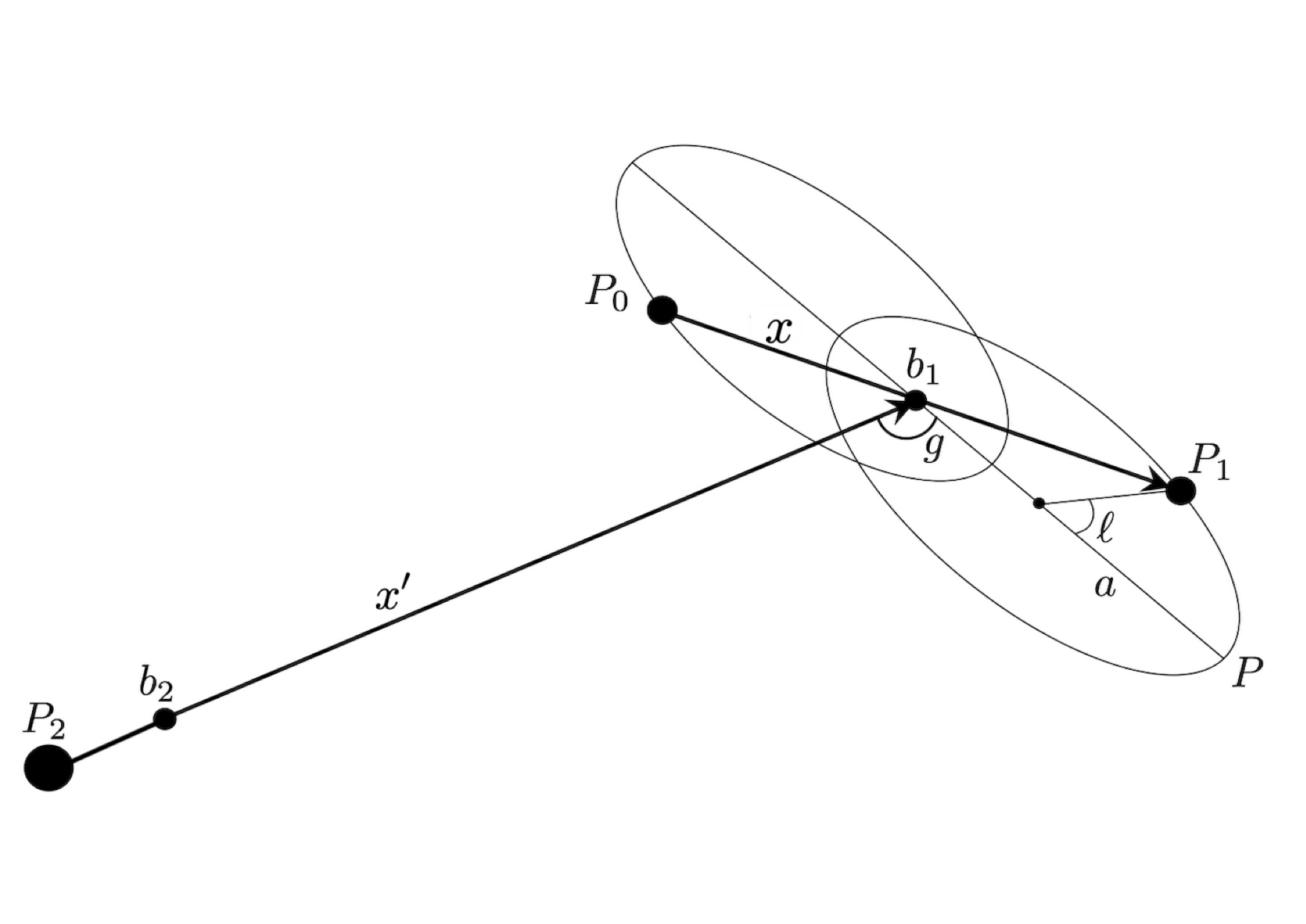}
\caption{\label{fig:model}Schematic representation of the model we are dealing with. The model is composed by three bodies $P_{0}, P_{1}, P_{2}$, where the first two have equal masses $m_{0}$ and the third body has the largest mass $\kappa  m_{0}$, $\kappa >1$. The point $b_{1}$ is the barycentre of $P_{0}$ and $P_{1}$, whilst $b_{2}$ is the barycentre of all the three points (close but different from $P_{2}$). 
}
\end{figure}

 The Hamiltonian governing the motions of three point masses undergoing Newtonian attraction is, as well known,
\begin{eqnarray}\label{first_Ham}
	\HH = \frac{{\vert y_0\vert}^2}{2 m_{0}} +
	\frac{{\vert y_1\vert}^2}{2 \mu m_0}+\frac{{\vert y_2\vert}^2}{2\kappa m_{0}}-\frac{\mu m_{0}^2}{\vert x_0-x_1\vert}-\frac{\kappa  m_{0}^2}{\vert x_0-x_2\vert}-\frac{\mu \kappa  m_{0}^2}{\vert x_1-x_2\vert}.
\end{eqnarray}
Here,  $x_{0}, x_{1}, x_{2}$ and $y_{0}, y_{1},y_{2}$ are, respectively, positions and impulses of the three particles relatively to a prefixed orthonormal frame $(i, j, k)\subset \R^3$; $m_0$, $m_1=\mu m_0$, $m_2=\kappa  m_0$, with $y_i=m_i\dot x_i$, are their respective gravitational  masses; $\vert\cdot\vert$ denotes the Euclidean distance and the gravity constant has been taken equal to one, by a proper choice of the unit system. In the sequel, in accordance to our problem, we shall take $x_i$, $y_i\in \R^2\times \{0\}\simeq\R^2$ and $\mu=1\ll \kappa $, so that $x_0$, $x_1$ correspond to the position coordinates of the asteroids; $x_2$ is the planet. The Hamiltonian $\HH$ is translation invariant, so we rapidly switch to a translation--free Hamiltonian by applying the well known Jacobi reduction. We recall that this reduction consists of using, as position coordinates, the centre of mass $r_0$ of the system (which moves linearly in time); the relative distance $x$ of two of the three particles; the distance $x'$ of the third particle with respect to the centre of mass of the former two. Namely,
\begin{eqnarray}\label{x0}
r_0=(x_0+\mu x_1+\kappa  x_2)(1+\mu+\kappa )^{-1}\ ,\ 
x=x_1-x_0\ ,\ x'=x_2-(x_0+\mu x_1)(1+\mu)^{-1}\ .
    \end{eqnarray} 
Note that, under the choice of the masses specified above, we are choosing  the asteroidal coordinate $x_0$ as the ``starting point'' of the reduction. This reverses a bit the usual practice, as $x_0$ is most often chosen as the coordinate of the most massive body; see Figure~\ref{fig:model}. 
At this point, the procedure is classical: the new impulses $(p_0$, $y$, $y'$) are uniquely defined by the  constraint of symplecticity, with $p_0$ (``total linear momentum'') being proportional to the velocity of the barycentre.  Choosing (as it is possible to do) a reference frame centred at, and moving with, $r_0$, so to have $r_0\equiv0\equiv p_0$, after a suitable rescaling, one obtains (see Appendix~\ref{Jacobi})
\begin{eqnarray} \label{ham_rescaled}
	\HH = \frac{1}{2 m_{0}} {\vert y\vert}^2 
	+
	\frac{\sigma }{2 m_{0}} {\vert y'\vert}^2 
	-\frac{m_{0}^2}{\vert x\vert} - \frac{m_{0}^2  \sigma}{\vert x' + \bar{\beta} x\vert}- \frac{\bar{\beta}}{\beta} \frac{m_{0}^2  \sigma }{\vert x' - \beta x \vert}, 
\end{eqnarray}
with
\begin{eqnarray}\label{beta}
\beta = \frac{\kappa ^2 (1+\mu)}{\mu^2 (1+\mu+\kappa )}\ ,\quad \bar{\beta} =\frac{\kappa ^2 (1+\mu)}{\mu (1+\mu+\kappa )}\ , \quad \sigma = \frac{\kappa ^3 {(1+\mu)}^2}{\mu^2 (1+\mu+\kappa )}.
\end{eqnarray}
The choice $\kappa \gg\mu =1$ gives $\beta=\bar \beta\gg1$ and simplifies $\HH$ to
\begin{eqnarray} \label{ham_mu_eq_1}
	\HH = \frac{ {\vert y\vert}^2}{2 m_{0}}  -\frac{m_{0}^2}{\vert x\vert}  
	+
	\frac{\sigma {\vert y'\vert}^2}{2m_{0}} - \frac{m_{0}^2  \sigma }{\vert x' + \beta x \vert}- \frac{m_{0}^2  \sigma }{\vert x' - \beta  x \vert}. 
\end{eqnarray}
\noindent From now on, we regard $\beta$ as  mass parameter, with $\beta\sim \kappa$ and $\sigma\sim \beta^2$. By choosing a region of the phase--space where 
\begin{eqnarray} \label{ph_sp_conditionOLD}
	\vert x'\vert>\vert \beta x\vert,
\end{eqnarray}
we ensure  the denominators of the two last terms in~\eqref{ham_mu_eq_1} to be different from zero. The Hamiltonian (\ref{ham_mu_eq_1}) with $x$, $x'$, $y$, $y'\in \R^2$, has $4$--degrees--of--freedom (\DOF \!, from now on), but is SO(2)--invariant. We choose a system of canonical coordinates which reduces this symmetry and hence lowers the number of \DOF to $3$. If $k=i\times j$ is normal to the plane of the orbits, we denote as
\begin{eqnarray*}
	C = x \times y \cdot k + x' \times y'\cdot k
\end{eqnarray*}
the total angular momentum, which is a constant of  the motion. Then we take a 3--\DOF system of coordinates, which we name $(\Lambda,G,R,\ell,g,r) \in \R^3\times{\mathbb T^2}\times\R$, where   $(\Lambda, G, \ell, g)$ are  ``Delaunay coordinates for the asteroid, relatively to $x'$'', while $(R, r)$ are  ``radial  coordinates  for the planet''.
More precisely, they are defined as
\begin{eqnarray}\label{coord}
\left\{
\begin{array}{lll}
 \Lambda  =  \sqrt{m_{0}^3  a} \\ 
 G = x \times y \cdot  k \\ 
 \ell  =   2 \pi \frac{\SC}{\SC_{\textrm{tot}}} \\ 
 g =  \alpha_{x',P}
\end{array}
\right. \, ,\qquad \left\{
\begin{array}{lll}
 R  =  y' \cdot \frac{x'}{\vert x'\vert}\\
 r   =  \vert x'\vert
\end{array}
\right.
\end{eqnarray}
where, considering the instantaneous ellipse generated by the first two terms in the Hamiltonian~\eqref{ham_mu_eq_1}, $a$ is the semi--major axis (see again Figure~\ref{fig:model}), $\SC$ and $\SC_{\textrm{tot}}$ are the area of the ellipse spanned from the perihelion $P$ and the total area and $\alpha_{x',P}$ is the angle between the direction of $x'$ and $P$ relatively to the positive direction established by $x\times y$. With these notations, $\ell$ represents the mean anomaly, $G$ is the projection of the angular momentum of the asteroid on the direction of the unit vector $k$ and $g$ is the anomaly of the perihelion $P$ with respect to the direction of $x'$. Using the coordinates (\ref{coord}), condition (\ref{ph_sp_conditionOLD}) becomes
\begin{eqnarray} \label{ph_sp_condition}
	\varepsilon<\frac{1}{2}\ ,\qquad \varepsilon:=\frac{\beta a}{r}
	\end{eqnarray}
 as a body moving on an ellipse does not go further than twice the semi--axis from the focus of the ellipse. The canonical character of the coordinates (\ref{coord})  has been discussed, in a more general setting, in~\cite{gPi19}. In terms of the coordinates (\ref{coord}), the  Hamiltonian~\eqref{ham_mu_eq_1} reads
\begin{eqnarray}
\label{eq:H}
	\HH &=&-\frac{m_{0}^5}{2\Lambda^2}
	+
	\frac{\sigma}{2m_{0}} \Big(R^2+ \frac{(G-C)^2}{r^2}\Big)
	-
	\frac{\sigma m_{0}^2}{\sqrt{r^2+2\beta a r p + \beta^2 a^2 \varrho^2}}
	\nonumber\\
	&&-
	\frac{\sigma m_{0}^2}{\sqrt{r^2-2\beta a r p + \beta^2 a^2 \varrho^2}} \, ,
\end{eqnarray}
where, for short, we have let
\begin{eqnarray*}
\varrho=\varrho(\Lambda,G,\ell)=1-e \cos \xi(\lambda, G, \ell)\ ,\quad 	 p = p(\Lambda,G,\ell,g)= (\cos \xi -e) \cos g
 -
 \frac{G}{\Lambda} \sin \xi \sin g.
\end{eqnarray*}
Here, 
\begin{eqnarray*}
	e=e(\Lambda,G)=\sqrt{1-\frac{G^2}{\Lambda^2}}
\end{eqnarray*}
is the eccentricity, and $\xi=\xi(\Lambda, G, \ell)$ denotes the eccentric anomaly, defined as the solution of  Kepler's equation 
\begin{eqnarray*}
  \xi - e(\Lambda, G) \sin \xi=\ell \, .
\end{eqnarray*}
 The next step is to switch to the 2--\DOF $\ell$--averaged (hereafter, {\it secular}) Hamiltonian, which we write as
\begin{eqnarray}\label{eq:AvH}
	\bar{\mathcal{H}}(G, R, g, r)&=&\frac{1}{2\pi} \int_{\mathbb{T}} \mathcal{H} \, \dd \ell\nonumber\\
	&=&-\frac{m_0^5}{2\Lambda^2}+
	\sigma \KK(R, r, G)
	+\sigma \UU(G, g, r)\, ,
\end{eqnarray}
with 
\begin{eqnarray}\label{KU}
\KK(R, r, G)&  := & \frac{R^2}{2m_0}+ \frac{(G-C)^2}{2m_0 r^2}-\frac{2  m_0^2}{r}\nonumber\\
 \UU(G, g, r)&:=&\UU_+(G, g, r)+\UU_-(G, g, r)+\frac{2 m_0^2}{r}
\end{eqnarray}
where
\begin{eqnarray}\label{averaged}
\UU_\pm (G, g, r):=-\frac{m_0^2}{2\pi}\int_0^{2\pi}\frac{d\ell}{\sqrt{r^2\pm2\beta a r p + \beta^2 a^2 \varrho^2}}.
\end{eqnarray}
In (\ref{eq:AvH}) we have omitted to write $\Lambda$ and $C$ among the arguments of $\bar\HH$, as they now play the r\^ole of parameters. Observe that the function $\UU$ is $\pi$--periodic in $g$, as changing $g$ with $g+\pi$ corresponds to swap $\UU_+$ and $\UU_-$, as one readily sees from (\ref{eq:H})--(\ref{averaged}). \\
We do not provide rigorous bounds ensuring that the secular problem may be regarded as a good model for the full problem. Heuristically, we expect that this is true  as soon as~\eqref{ph_sp_condition} is strengthened requiring, also,
\begin{eqnarray}\label{small pert}r\gg \beta^{3/2} a.
\end{eqnarray}
Indeed, extracting $r$ from the denominators of the two latter functions in (\ref{eq:H})  and expanding the resulting functions in powers of $\frac{\beta a}{r}$, one sees that the lowest order terms depending on $\ell$ have size $\frac{m_0^2\sigma\beta a}{r^2}\sim {m_0^2}\frac{\beta^3 a}{r^2}$ (recall that $\sigma\sim \beta^2$). So, such terms are negligible compared to the size $\frac{m_0^2}{a}$ of the Keplerian term, provided that~\eqref{small pert} is verified.\\
Neglecting the constant term $-\frac{m_0^5}{2\Lambda^2}$  and, after a further change of time, the common factor $\sigma$ in the remaining terms, the secular Hamiltonian (\ref{eq:AvH}) reduces to
\begin{eqnarray}\label{averagedHamNEW}
\hat\HH(G, R, g, r)=\KK(R, r, G)+\UU(G, g, r)\, .
\end{eqnarray}
We now specify the range of parameters  $C$, $\Lambda$ and $\beta$ and the region of the phase space for the coordinates $(G, R, g, r)$ that we consider in this paper. In particular, we look for values of parameters and coordinates where the   Hamiltonian (\ref{averagedHamNEW}) is {\it weakly coupled}, and describe the motions we expect to find in such region. As above, our discussion will be extremely informal.\\
First of all, we  take  $\Lambda$ and $C$ verifying
\begin{eqnarray}\label{C and Lambda}
\Lambda\ll C\, .
\end{eqnarray}
This condition implies that also $|G|\ll C$ (as $|G|<\Lambda$) and hence $\KK$ affords the natural splitting $\KK=\KK_0+\KK_1$, where
  \begin{eqnarray*}
\KK_0=\frac{R^2}{2 m_0}+\frac{C^2}{2 m_0  r^2}-\frac{2 m_0^2}{r}\ ,\qquad \KK_1=\frac{G(G-2C)}{2 m_0 r^2} \, .
\end{eqnarray*}
 We consider a region of phase--space where $r$ and $R$ take values 
  \begin{eqnarray}\label{circular}r\sim r_0=\frac{C^2}{2m_0^3}\ ,\quad R\sim 0 \, .\ 
  \end{eqnarray}
 These are the values where  $\KK_0$ attains its minimum, and correspond to circular motions of the planet, with $r_0$ being the radius of the circle.  In the region of phase space defined by (\ref{circular}),  the relative sizes of $\KK_1$ and $\UU$ to $\KK_0$ are
 \begin{eqnarray}\label{weights}
 \|\KK_1\|< c_1\,\frac{\Lambda}{C} \|\KK_0\|\ ,\qquad  \|\,\UU\|<  c_2\,\varepsilon^2 \|\KK_0\| \, ,
 \end{eqnarray}
where $c_i$ are independent of $m_0$, $\beta$, $\Lambda$ and $C$. Even though (by (\ref{C and Lambda}) and (\ref{small pert})) $\KK_1$ and $\UU$ are small compared to $\KK_0$, however, they cannot be neglected, as their sum  governs the slow motions of the coordinates $G$ and $g$, which do not appear in $\KK_0$. Remark that $\KK_1$ and $\UU$ are coupled with $\KK_0$, since they depend on $r$. It is however reasonable to expect that, as long as the minimum of $\KK_0$ cages $r$ to be close  to the value $r_0$, the coupling is weak and the dynamics of  $G$ and $g$ is, at a first approximation, governed by the 1 \DOF Hamiltonian \begin{eqnarray}\label{F}\FF(G, g):=(\KK_1+\UU)|_{r= r_0}\ .\end{eqnarray} To understand the global phase portrait of $\FF$ in the plane $(g, G)$, we need to recall some results from~\cite{gPi20b}. We go back to the functions $\UU_\pm$ in (\ref{averaged}), which enter in the definition of $\UU$. In~\cite[Section 3]{gPi20b}, it is proved that, under the assumption (\ref{ph_sp_condition}), the following identity holds
\begin{eqnarray}\label{identity}\UU_\pm(G, g, r)=-\frac{m_0^2}{2\pi r}\int_0^{2\pi}\frac{(1-\cos\xi)d\xi}{\sqrt{1\mp \varepsilon (1-\cos\xi) t_\pm+\varepsilon^2 (1-\cos\xi)^2}}
\end{eqnarray}
 with $\varepsilon$ as in (\ref{ph_sp_condition}) and
 \begin{eqnarray*} t_\pm(G, g, \varepsilon):=\sqrt{1-\frac{G^2}{\Lambda^2}}\cos g\pm \varepsilon \frac{G^2}{\Lambda^2}\ . 
\end{eqnarray*}
 \begin{figure}
 \includegraphics[height=5.0cm, width=6.0cm]{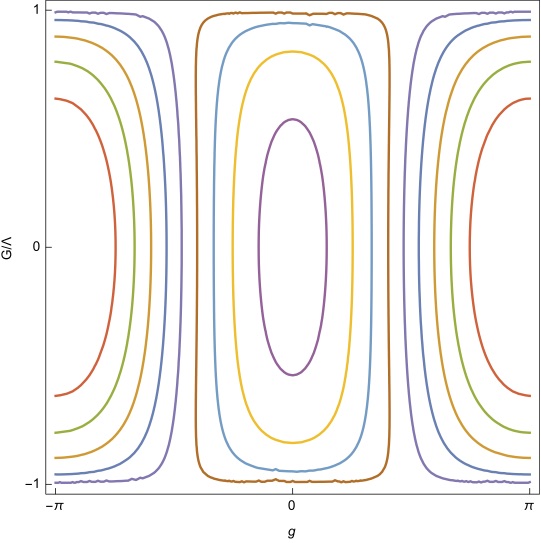} 
 \includegraphics[height=5.0cm, width=6.0cm]{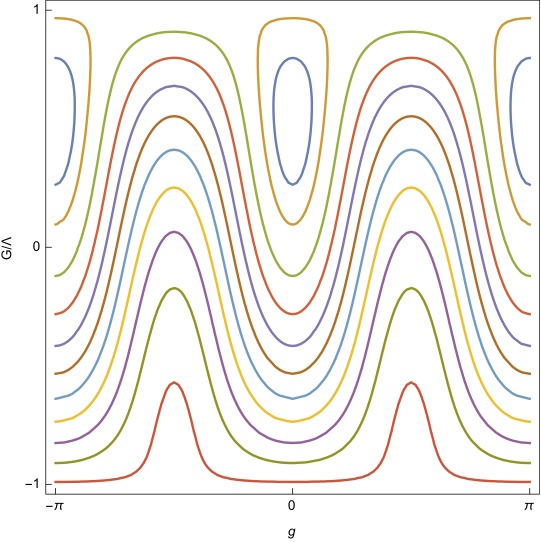}
 \caption{Left: the phase portrait of $t_+(\cdot, \cdot, \varepsilon)$ in the plane $(g, G/\Lambda)$, for $0<\varepsilon<\frac{1}{2}$. Right: the  phase portrait of $\FF$ in the plane $(g, G/\Lambda)$, with $m_0$, $C$, $\Lambda$ and $\beta$ as in (\ref{parameters}).}\label{unperturbed phase portrait}
 \end{figure}
 The equality (\ref{identity}) has two main consequences. The former is that, even though the transformation 
 (\ref{coord}) looses its meaning when $G=0$, however, $\UU_\pm$ keep their regularity, provided that (\ref{ph_sp_condition}) holds. Indeed, the functions $t_\pm$ are regular at $G=0$ and, being bounded below by $-1$  and  above  by $1$,  the denominator of the function under the integral never vanishes, under (\ref{ph_sp_condition}), as it is immediate to verify.   
Secondly, the phase portrait of the functions $\UU_+(\cdot, \cdot, r)$, $\UU_-(\cdot, \cdot, r)$  coincides, a part for a rescaling,  with the one of $t_+(\cdot, \cdot, \varepsilon)$, $t_-(\cdot, \cdot, \varepsilon)$, respectively. In particular,
 $\UU_+(\cdot, \cdot, r)$ and $\UU_-(\cdot, \cdot, r)$ have elliptic equilibria at $(G, g)=(0, 0)$ and $(G, g)=(0, \pi)$, because this is true for $t_\pm$, as it is immediate to check.
The phase portrait of $t_+(\cdot, \cdot, \varepsilon)$ for $\varepsilon<\frac{1}{2}$ is  shown in Figure~\ref{unperturbed phase portrait} (left); the one of $t_-(\cdot, \cdot, \varepsilon)$ is  specular, interchanging the equilibria. We now merge these informations, in order to build the phase portrait of the function $\FF$ in (\ref{F}).
By the Implicit Function Theorem, one can argue that,  for an open set of values of the parameters, due to the linear term in $G$ in $\KK_1$, the equilibria of $\UU_+$ and $\UU_-$ are shifted along the $G$--axis, but are not  destroyed. Quantifying this shift is not easy, as  $\UU$ has an involved dependence on $t_+$, $t_-$. Based on the $\varepsilon$--expansion of $\UU$, with
   \begin{eqnarray}\label{parameters}
m_0=1\ ,\quad C= 75\ ,\quad \Lambda=\sqrt{a}= 3\ ,\quad \beta= 40 
   \end{eqnarray}
   (which comply with~\eqref{ph_sp_condition},~\eqref{small pert},~\eqref{C and Lambda})
we obtain the 
 phase portrait of $\FF$ as in Figure~\ref{unperturbed phase portrait} (right). We observe that, at contrast with the figure at left--hand side, where the motions are purely of elliptic kind, the phase portrait at right--hand side also includes rotational motions. The linear term of $\KK_1$ is responsible of this fact, breaking the symmetry $G\to -G$. We underline at this respect that  the present framework  is in a sense complementary to the one studied in~\cite{gPi20b}, where   the phase portrait of $\FF$  has, in fact, only elliptic motions: in that case, the linear term of $\KK_1$ does not exist, as $C$ is fixed to $0$. Remark also that the vanishing of $C$ in~\cite{gPi20b} affects
condition (\ref{C and Lambda}) (which is not satisfied)  and the motions generated by $\KK_0$ (which are collisional, rather than circular). 
\\
The purpose of this paper is to show that, if the parameters are chosen about  (\ref{parameters}) and the energy is  fixed to the level of a suitable initial datum $(G_\star, R_\star, g_\star, r_\star)$ satisfying (\ref{circular})  (see Appendix~\ref{Choice of the parameters} for the exact values), then, in the system  (\ref{averagedHamNEW}) a topological horseshoe wakes up in the plane $(g, G/\Lambda)$. The analysis will be purely  numerical, based on techniques developed in~\cite{aGi19},~\cite{ZGLICZYNSKI200432}. More details on the methodological strategy are given along the following sections.


\section{Poincar\'e mapping} \label{sec:Poi_map}
\begin{figure}
\includegraphics[height=8.0cm, width=8.0cm]{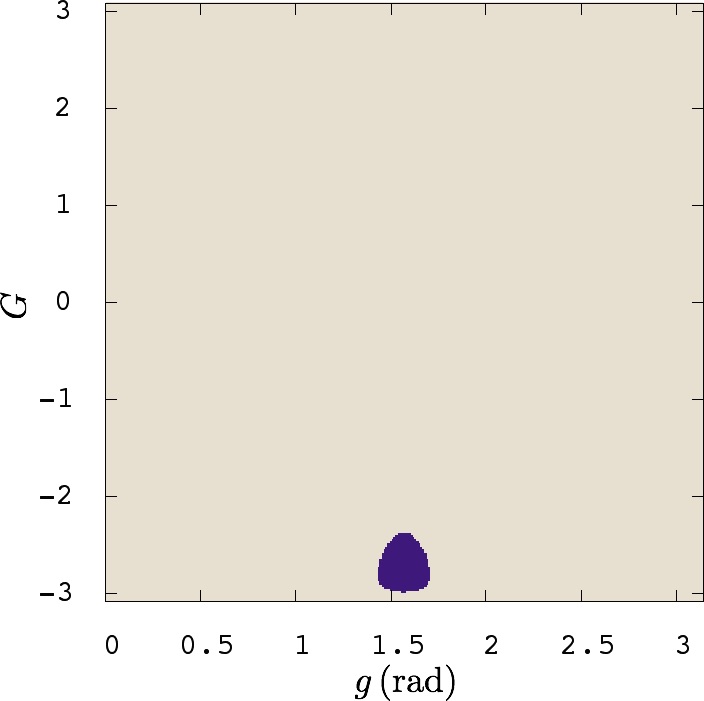}
\caption{\label{fig:fig1J} The admissible points of the $(g,G)$ section are displayed in cream colour. They correspond to points satisfying the energetic condition $\mathcal{H}=h_{\star}$ with $R^2 \ge 0$.
The complementary set (points leading to negative $R^2$) appear in purple and define the \textit{inadmissible seeds}. See text for more details.}
\end{figure}

\begin{figure}
\includegraphics[height=13.0cm, width=13.0cm]{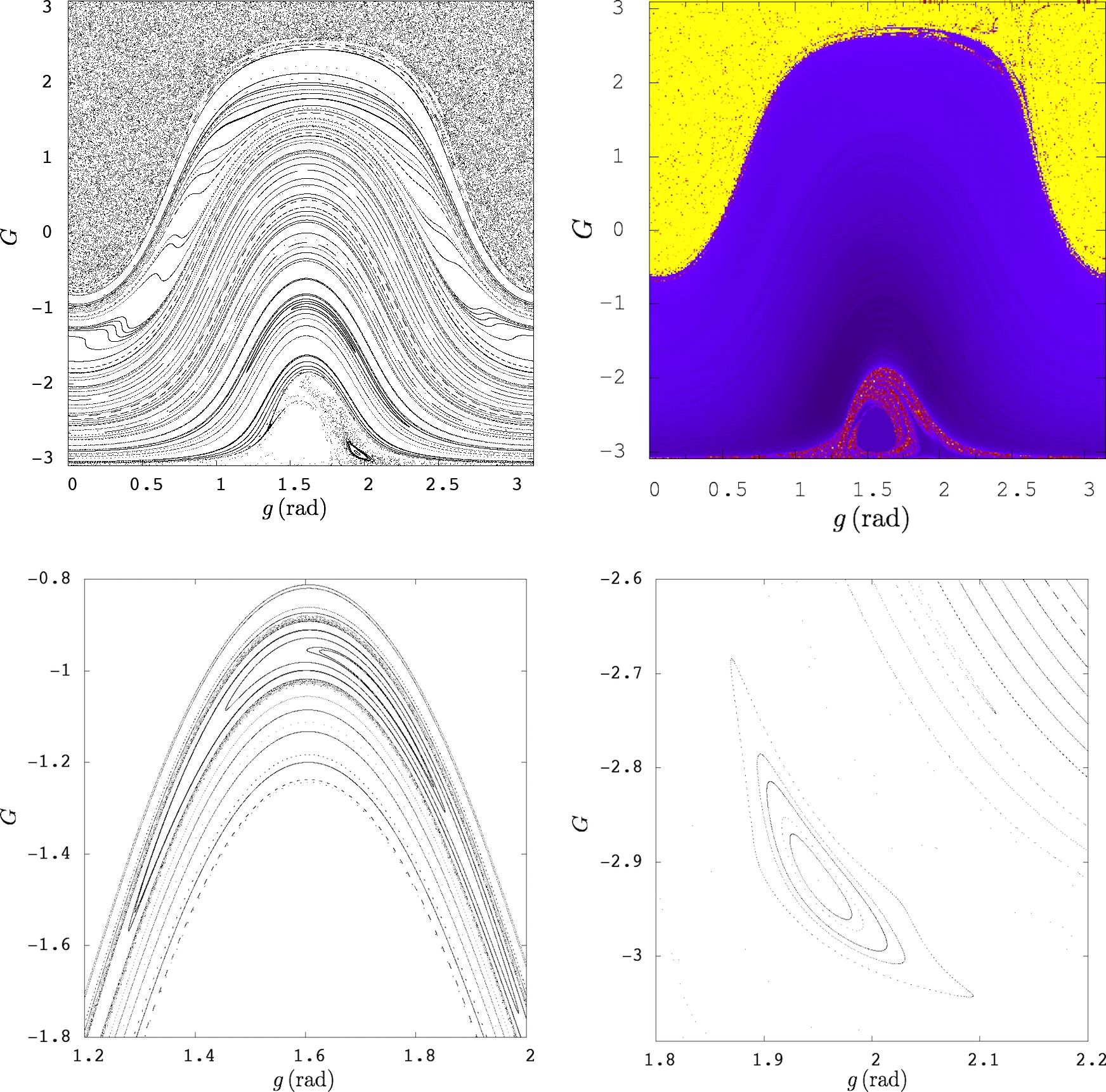}
\caption{\label{fig:fig2J} Phase--space structures of the mapping 
$P$ at different scales. 
Upper left: global phase--space;
lower: microscales structures;
upper right: the global phase--space analysis obtained by iterating the mapping $P$ is confirmed by computing finite time chaos indicators based on the variational dynamics derived from the continuous model $\mathcal{H}$.
}
\end{figure}

\noindent From now on, we neglect to write the ``hat'' in (\ref{averagedHamNEW}). Moreover, for the purposes of the computation, we replace the function $\UU$   with a finite sum 
\begin{eqnarray}\label{Uk}
\UU_k=\sum_{\nu=1}^{k} q_{\nu}(G,g,r) \Big(\beta \frac{a}{r}\Big)^{\nu}
\end{eqnarray}
where  $q_{\nu}(G,g,r)$ are the Taylor coefficients in the expansion of $\UU$ with  $\nu=1$, $\ldots$, $k$. Using the parity of $\UU$ as a function of $r$, these coefficients have the form
\begin{eqnarray*}
q_{\nu}(G,g,r) =\left\{
\begin{aligned}
	&\frac{m_0^2}{r} \sum_{p=0}^{\nu/2} \tilde q_{p}(G) \cos (2p\,g) \quad&{\rm if}\ \nu\ {\rm is\ even} \,   
	\\
	 &\ 0 \quad &{\rm otherwise}\ .
\end{aligned}
\right.	
\end{eqnarray*}
In our numerical implementation, we use the truncation in~\eqref{Uk} with  
$k=k_{\max}=10$, so as to balance  accuracy and number of produced terms. 
We still denote as $\HH$ the resulting Hamiltonian:
\begin{eqnarray} \label{eq:H_fin} \begin{aligned} 
\HH (G,R,g,r)   & =   \KK(G,R,r) +\UU_k(G,g,r) \\ 
			& = \frac{1}{2m_{0}} \Big(R^2+ \frac{(G-C)^2}{r^2}\Big)-\frac{2m_0^2}{r}
+
\sum_{\nu=1}^{k} q_{\nu}(G,g,r) \Big(\beta \frac{a}{r}\Big)^{\nu}. \,
\end{aligned}
\end{eqnarray}
The study of the  secular $2$--\DOF Hamiltonian in the continuous time $t$ can be reduced to the study of a discrete mapping through the introduction of ad--hoc Poincar\'e's section~\cite{jMe92}. The advantage consists in reducing further the dimensionality of the phase--space, and, in the case of $n=2$, to sharpen the visualisation of the dynamical system. In fact, for a $2$--\DOF system, the phase--space has dimension $4$ and, due to the conservation of the energy (the Hamiltonian $\mathcal{H}$ itself), orbits evolve on a three--dimensional manifold $M$. By choosing an appropriate surface $\Sigma$ transverse to the flow, one can look at the intersections of the orbits on the intersection of $M \cap \Sigma$, \ie a two--dimensional surface. The surface $\Sigma$ chosen is a plan passing through a given point $(G_{\star},g_{\star},r_{\star}) $ and normal to the associated orbit, \ie to the velocity vector $(v_{G}^{\star},v_{g}^{\star},v_{r}^{\star})$; it is defined by
\begin{eqnarray*} 
	\Sigma
	=
	\big\{(G,g,r): v_{G}^{\star} (G-G_{\star})+ v_{g}^{\star} (g-g_{\star})+ v_{r}^{\star} (r-r_{\star}) = 0\big\}.
\end{eqnarray*}
Let us now formally introduce the Poincar\'e map. We start by defining two operators $l$ and $\pi$ consisting in ``lifting'' the initial two--dimensional seed $z=(G,g)$ to the four--dimensional space $(G,R,g,r)$ and ``projecting'' it back to plan after the action of the flow--map $\Phi^{t}_{\mathcal{H}}$ during the first return time $\tau$. The lift operator reconstructs the four--dimensional state vector from a seed on $D \times \mathbb{T}/2$, where the domain $D$ of the variable $G$ is a compact subset of the form $[-\Lambda,\Lambda]$. For a suitable $(\mathcal{A},A) \subset \R^{2} \times \R^{2} $, its definition reads 
\begin{eqnarray}
l: \quad  D \times \mathbb{T}/2 \supset \mathcal{A}   & \rightarrow & D  \times \mathbb{T}/2 \times A \notag \\
   	      z	     & \mapsto &  \tilde{z}=l(z) \notag,
\end{eqnarray}
where  $\tilde{z}=(G,g,R,r)$ satisfies the two following conditions:
\begin{enumerate}
	\item \textit{Planarity condition.} The triplet $(G,g,r)$ belongs to the plane $\Sigma$, \ie $r$ solves the algebraic condition 
	$v_{G}^{\star} (G-G_{\star})+ v_{g}^{\star} (g-g_{\star})+ v_{r}^{\star} (r-r_{\star}) = 0$.  
	\item \textit{Energetic condition.} The component $R$ solves the energetic condition $\mathcal{H}(G_{\star},g_{\star},R_{\star},r_{\star})=h_{\star}$. The Hamiltonian is separable in $R$, so this condition amounts to solve a quadratic equation. If $R^2 \ge 0$, then we choose the root associated to the ``positive'' branch $+\sqrt{R^2}$. If $R^2 <0$, then we are led to the notion of \textit{inadmissible seed}. The set of admissible seeds, noted by $\mathcal{A}$, for the chosen section $\Sigma$ is portrayed in Figure\,\ref{fig:fig1J}.  
\end{enumerate} 
The projector $\pi$ is the projection onto the first two components of the vector,
\begin{eqnarray*}
\pi: \quad  D \times \mathbb{T}/2 \times A &\rightarrow & D \times \mathbb{T}/2  \notag \\
   	        \tilde{z}=(z_{1},z_{2},z_{3},z_{4})	     &\mapsto  & \pi(z)=(z_{1},z_{2})\, .
\end{eqnarray*}
The Poincar\'e mapping is  therefore defined and constructed as 
\begin{eqnarray}
P: \quad D \times \mathbb{T}/2 &\rightarrow & D \times \mathbb{T}/2 \notag \\
   	      z	     &\mapsto &  z'=P(z)= 
	      \big(\pi \circ \Phi^{\tau(z)}_{\mathcal{H}} \circ l \big)
	      (z)\, .\notag
\end{eqnarray}
 The mapping is nothing else than a ``snapshots'' of the whole flow at specific return time $\tau$. It should be noted that the successive (first) return time is in general function of the current seed (initial condition or current state), \ie $\tau=\tau(\tilde{z})$, formally defined (if it exists) as 
\begin{eqnarray*}
	\tau(z) = \inf \Big\{ t \in \mathbb{R}_{+}, 
	\big(G(t),g(t),r(t)\big)	\in \Sigma   \Big\},
\end{eqnarray*}
where $\big(G(t),g(t),r(t)\big)$ is obtained though $\Phi^{t}_{\mathcal{H}}(G,R,g,r)$. The Poincar\'e return map we described has been constructed numerically based on the numerical integration of the Hamiltonian equation of motions (the details regarding our numerical settings are presented in the Appendix~\ref{App:NS}.)  This mapping being now explicit, we are able to unveil the phase--space structures through successive iterations of $P$. Figure~\ref{fig:fig2J} presents the successive coordinates of $\{P^{n}(z)\}$ where the initial seeds $z$ cover a discretisation of $D \times \mathbb{T}/2$ domain (mesh) and $n \sim10^3$. The phase--space structures can be roughly categorised in three distinct zones. In the lower part, say for $G < -2$, we can distinguish one ``pic'' centred around $g= \pi/2$. One elliptic zone is immersed inside this structure, surrounded by ``scattered dots'', indicative of chaos.  There is a large region of the phase--space foliated by circulational tori. The last upper region is a large zone where almost all regular structures have disappeared. The panel provided by Figure\,\ref{fig:fig2J} presents some magnifications of phase--space structures. The obtained phase--space structures have been confirmed using a finite time dynamical chaos indicator, the Fast Lyapunov Indicator (FLI) computed with the whole flow on an iso--energetic section (see Appendix~\ref{App:FLI} for more details). The FLIs computation relies on monitoring the growth over time of the \textit{tangent vector} under the action of the tangent flow--map (\textit{variational dynamics}). The final FLIs values are colour coded according to their values and projected onto the section to provide a stability chart. Stable orbits correspond to dark regions, orbits possessing the sensitivity to initial conditions appear in reddish/yellow color.  As shown in Figure\,\ref{fig:fig2J}, the FLIs confirm nicely the global structures depicted via the mapping. Moreover, numeric suggests that the lift of $P$ on the variables $(G, g, r)$ (\ie the map obtained from $\Phi_{\HH}^{\tau}$ by projection on $(G, g, r)$) is generically twist.

\subsection{Hyperbolic structures and heteroclinic intersections} 
Equilibrium points of the mapping $P$ (\ie periodic orbits of the Hamiltonian system (\ref{eq:H_fin})), have been found using a Newton algorithm with initial guesses distributed on a resolved grid of initial conditions in $D \times \mathbb{T}/2$ (again, see Appendix~\ref{App:NS} for more details regarding the numerical setup). We found more than $20$ fixed points $x_{\star}$ whose coordinates have been reported in Appendix~\ref{sub:CoordFP}. 
 The eigensystems associated to the fixed points have been computed to determine the local stability properties. The point $x_{\star}$ is hyperbolic when one of its real eigenvalues has modulus greater than one, the other less than one (expanding and contracting directions, respectively). In the case of complex eigenvalues, the point is elliptical. The result of the analysis is displayed on Figure 5 along with the following convention: hyperbolic fixed points appear as red crosses, elliptical points are marked with blue circles. As intuitively expected, the hyperbolic points are embedded within the chaotic sea. On the contrary, the stable islands host the elliptic points. Note that even the fixed--point in the small stability island has been recovered with the Newton scheme. In the vicinity of the unstable fixed--points, the dynamics is dominated by the stable and unstable manifolds who have the eigenvectors of $DP(x_{\star})$ asymptotically tangents near $x_{\star}$. The local stable manifold associated to an hyperbolic point $x_{\star}$,
\begin{eqnarray*}
	\mathcal{W}^{s}_{\textrm{loc.}}(x_{\star}) = 
	\Big\{ 
	x \, \vert \,
	\norm{P^{n}(x)-x_{\star}} \to 0, \, n \in \mathbb{N}_{+}, n \to \infty
	\Big\},
\end{eqnarray*}
can be grown by computing the images of a fundamental domain $I \subset E_{s}(x_{\star})$,  $E_{s}(x_{\star})$ being the stable eigenspace associated to the saddle point $x_{\star}$. We considered the simplest parametrisation of $I$, namely a normalised version of the eigenvector associated to the saddle point $x_{\star}$. This allowed us to compute a piece of $\mathcal{W}^{s}_{\textrm{loc.}}(x_{\star})$ under the action of the flow--map~\cite{cSi90,bKr06}. To compute the unstable manifold, the same computations are performed by reversing the time and changing $E_{s}$ by $E_{u}$. Finite pieces of those manifolds are presented in Figure~\ref{figJ:FPMan} for two saddle points. Following the well established conventions of the cardiovascular system (as reported in~\cite{jMe08}), the stable manifolds are displayed with blue tones, unstable manifolds appear in red tones. As we can observe, those curves intersect transversally forming the sets of heteroclinic points, trademark of the heteroclinic tangle and chaos~\cite{aMo02}.  We now have at hands all the necessary ingredients and tools to prove the existence of symbolic dynamics using covering relationships and their images under $P$. 

\begin{figure}
\includegraphics[height=9cm, width=10.0cm]{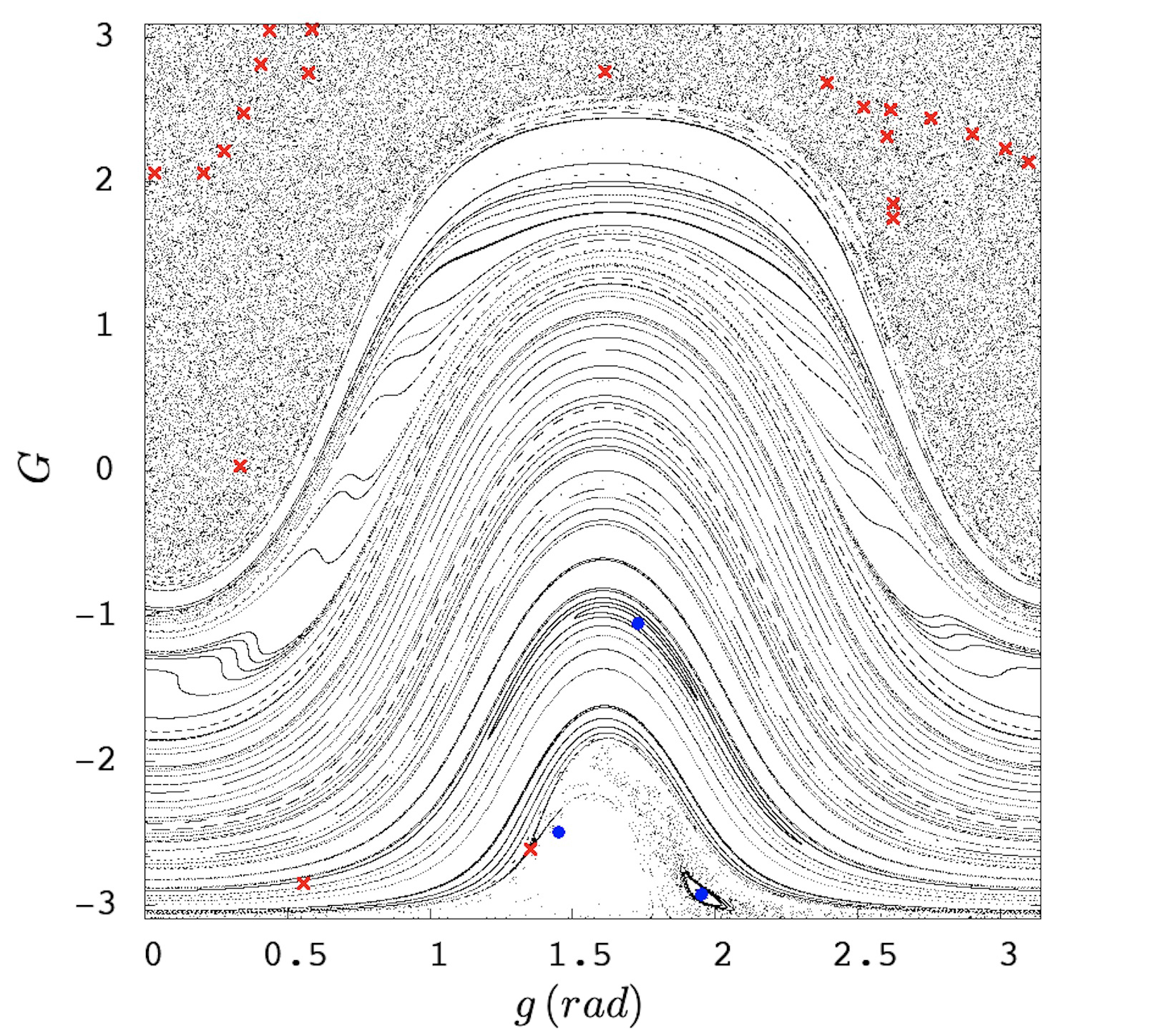}
\caption{\label{figJ:StabilityFP} 
Phase--space of $P$ together with its fixed--points. Hyperbolic points appear with red crosses, elliptical points appear with blue circles.
}
\end{figure}

\begin{figure}
\includegraphics[height=9cm, width=10.0cm]{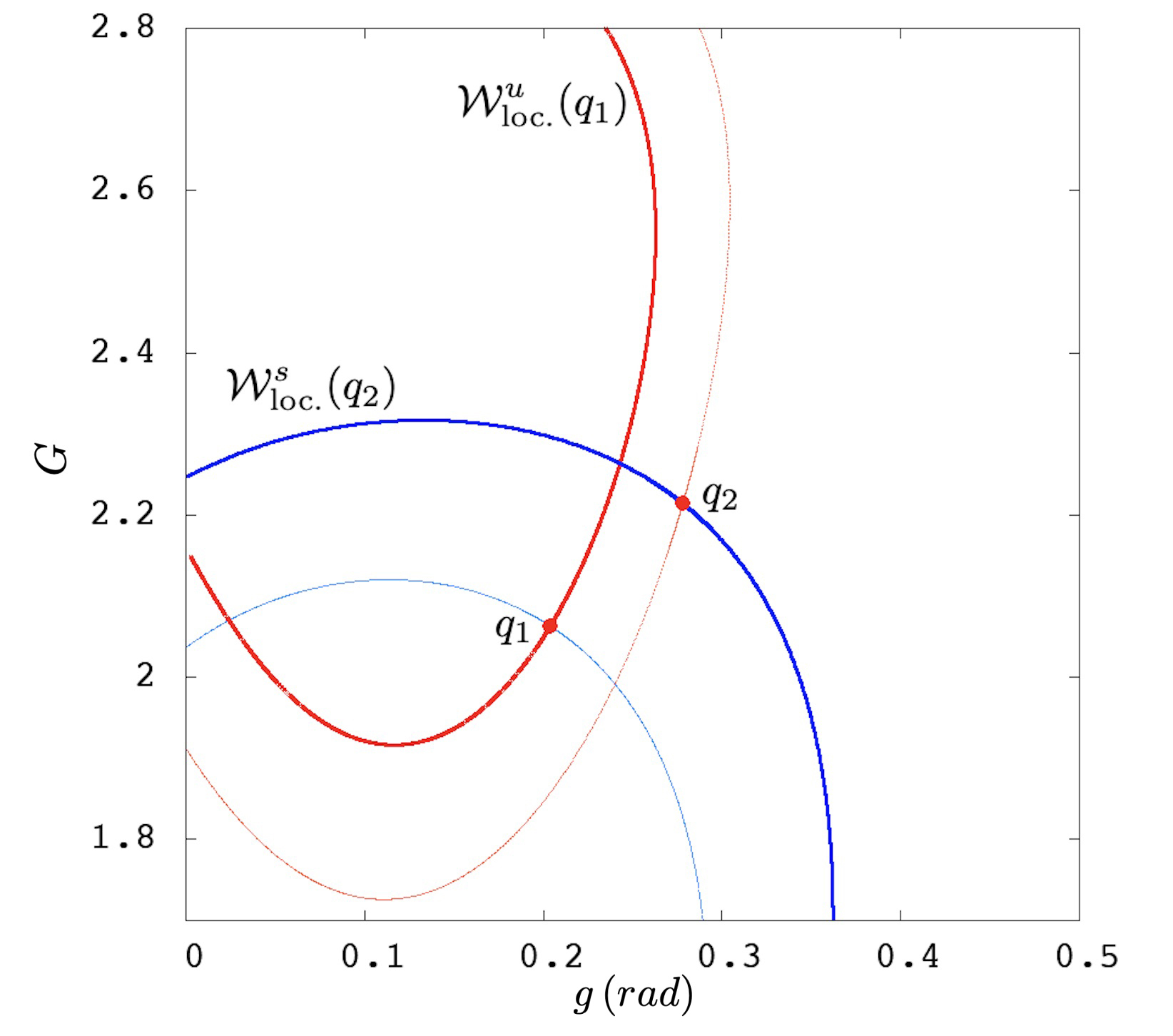}
\caption{\label{figJ:FPMan} 
Finite pieces of manifolds of two hyperbolic fixed points $q_1,q_2$. Their stable and unstable manifolds intersect transversally in (more than one) heteroclinic points. Stable manifolds are in blue while unstable manifolds are in red.
}
\end{figure}

\section{Symbolic dynamics via covering relations}
\noindent In this section we prove the existence of symbolic dynamics for the considered  model. The tools rely on ad--hoc covering relations that we present briefly following~\cite{aGi19}, in particular for the case $n=2$. 

\subsection{Covering relations and topological horseshoe}  
\begin{figure}
\includegraphics[height=3cm, width=6cm]{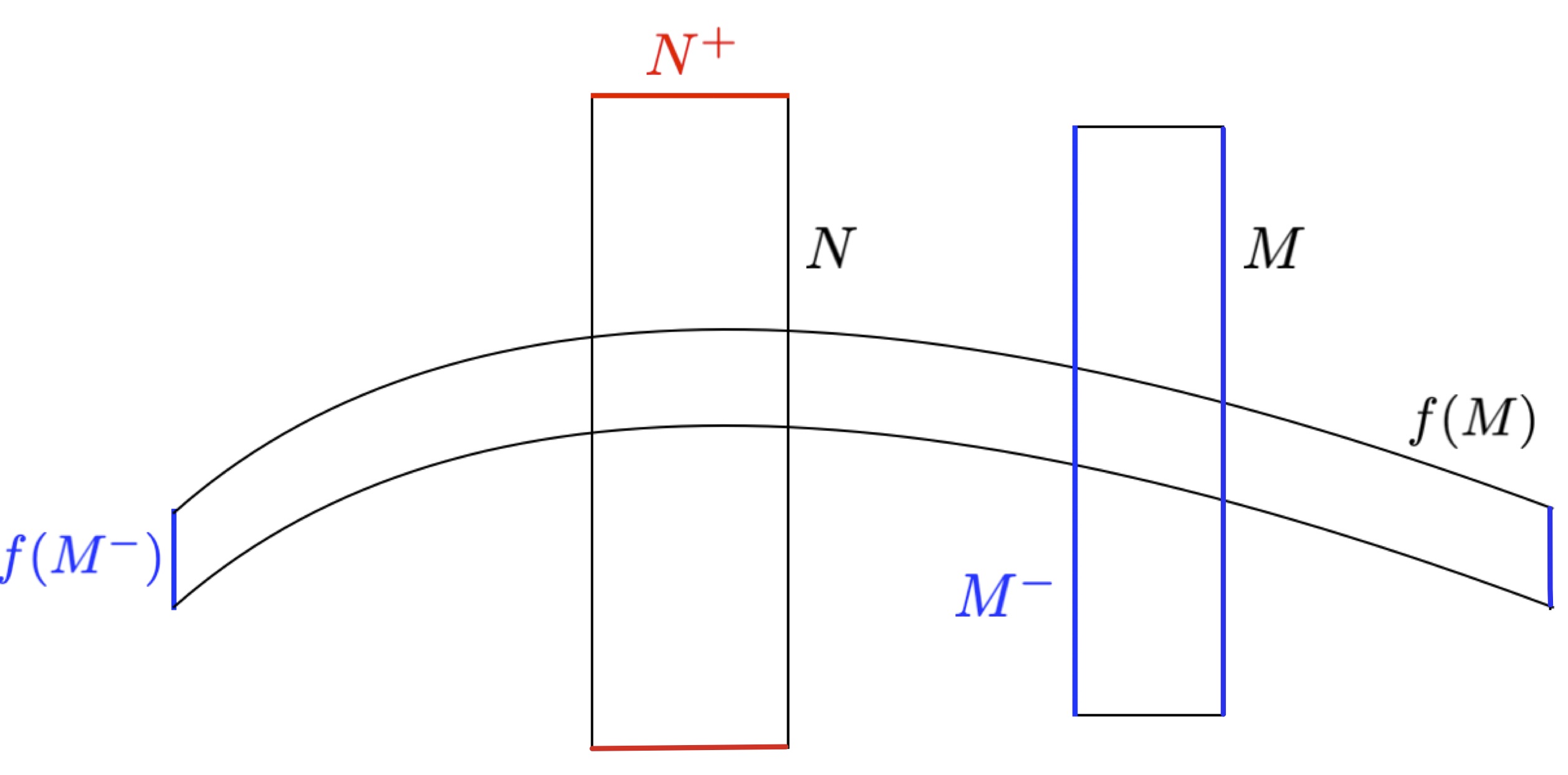} \qquad
\includegraphics[height=2.8cm, width=6cm]{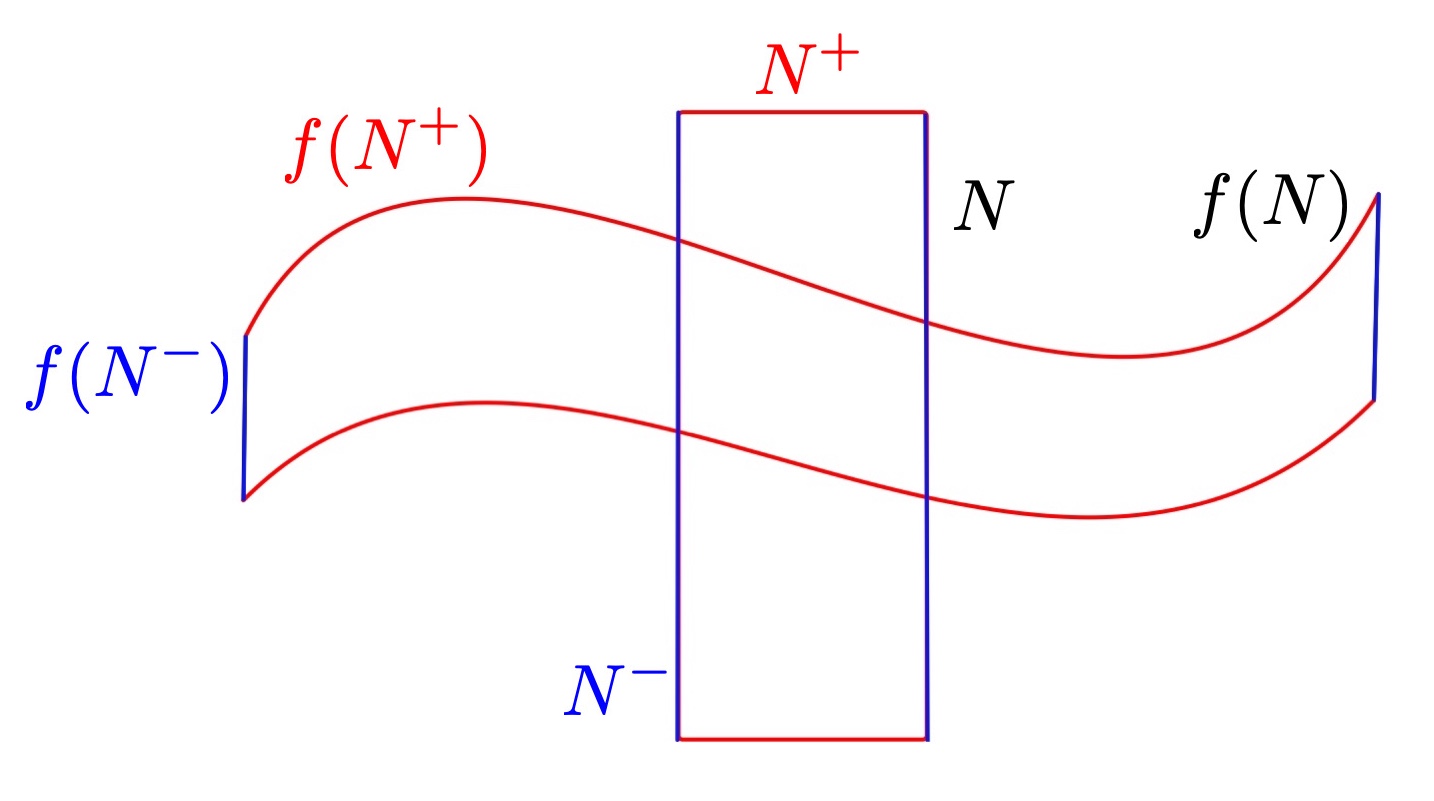}
\caption{\label{fig:covering}
Examples of covering relations. On the left $M \stackrel{f} \Longrightarrow N$. On the right, a case of self--covering $N \stackrel{f} \Longrightarrow N$ is illustrated. In red the entry sets and their image are represented, while in blue the exit sets and their images are represented.
}
\end{figure}

Let us introduce some notations. Let $N$ be a compact set contained in 
$\R^2$ and $u(N)=s(N)=1$ being, respectively, the {\it exit} and {\it entry dimension} (two real numbers such that their sum is equal to the dimension of the space containing $N$); let $c_{N} : \R^2 \rightarrow \R^2$ be an homeomorphism such that $c_N(N) = [-1,1]^2$; let $N_c=[-1,1]^2$, $N_c^-=\{-1,1\} \times [-1,1]$, $N_c^+= [-1,1]\times \{-1,1\}$; then, the two set $N^-= c_N^{-1}(N_c^-)$ and $N^+= c_N^{-1}(N_c^+)$ are, respectively, the {\it exit set} and the {\it entry set}. In the case of dimension $2$, they are topologically a sum of two disjoint intervals.  The quadruple $(N,u(N),s(N),c_N)$ is called a {\it h--set} and $N$ is called {\it support} of the $h$--set. Finally, let $S(N)_c^l =(-\infty, -1) \times \R$, \,  $S(N)_c^r =(1,\infty) \times \R$, and $S(N)^l = c_N^{-1}(S(N)_c^l) , \, S(N)^r = c_N^{-1}(S(N)_c^r) $ be, respectively, the left and the right side  of $N$.   The general definition of covering relation can be found in~\cite{aGi19}. Here we provide a simplified notion, suited to the case that $N$ is two--dimensional, based on\footnote{More precisely, Definition \ref{def_covering} is based on the {\it proof} of \cite[Theorem 16]{ZGLICZYNSKI200432}. Indeed,
 \cite[Theorem 16]{ZGLICZYNSKI200432} asserts that under conditions (1), (3) and one of the inclusions in \cite[(78) or (79)]{ZGLICZYNSKI200432}, one has $M \stackrel{f} \Longrightarrow N $ in the sense of~\cite{aGi19}. However, during the proof of 
\cite[Theorem 16]{ZGLICZYNSKI200432}, inclusions \cite[(78) or (79)]{ZGLICZYNSKI200432} are  only used to check the validity of (2).
} \cite[Theorem 16]{ZGLICZYNSKI200432}.
\begin{definition}
\label{def_covering}
	Let $f : \R^2 \rightarrow \R^2$ be a continuous map and $N$ and $M$ the supports of two $h$--sets. We say that $M$ $f$--covers $N$ and we denote it by $M \stackrel{f} \Longrightarrow N $ if:
	\begin{itemize}
	\item[(1)] $\exists\, q_0\in [-1, 1]$ such that $f(c_N([-1, 1]\times \{q_0\}))\subset {\rm int}( S(N)^l \bigcup N \bigcup S(N)^r)$,
	\item[(2)] $f(M^-) \bigcap N  = \emptyset$,
	\item[(3)] $f(M) \bigcap N^+  = \emptyset$.
\end{itemize}	
Conditions (2) and (3) are called, respectively, {\it exit} and {\it entry condition}.
\end{definition}
\noindent The case of self--covering is not excluded. The Figure\,\ref{fig:covering}  shows two schematic examples of covering relation between two different sets $N, M$ and a self--covering relation of $N$. The notions of covering relationships are useful in defining \textit{topological horseshoe} (confer~\cite{aGi19, ZGLICZYNSKI200432}).\\

\begin{figure}
\includegraphics[height=6cm, width=7cm]{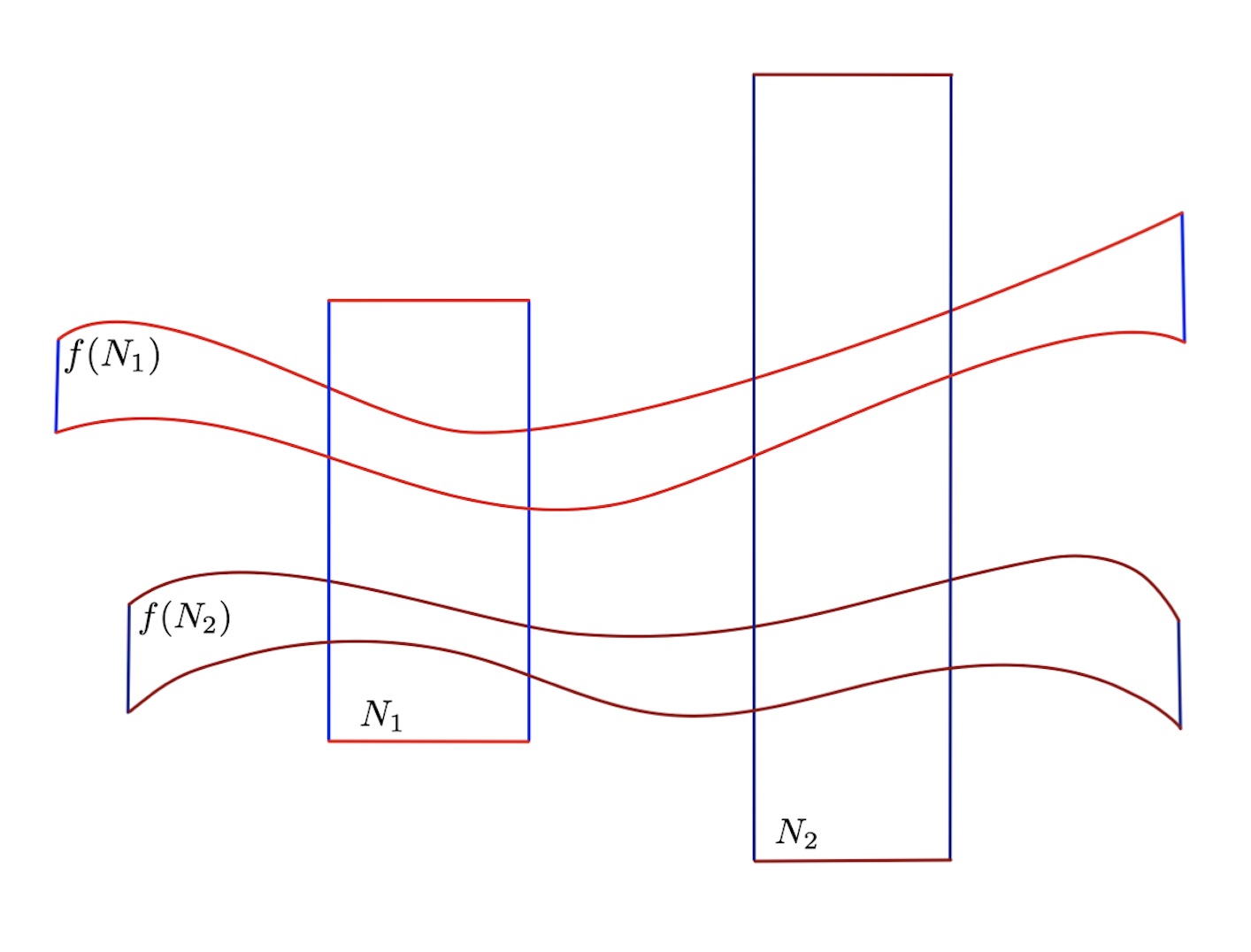}
\caption{\label{fig:horseshoe}
An example of topological horseshoe where both $N_{1},N_{2}$ cover themselves and each others. In red tones the entry sets and their images are represented, while in blue tones the exit sets and their images are represented.}
\end{figure}

\begin{definition}
	Let $N_{1}$ and $N_{2}$ be the supports of two disjoint $h$--sets in $\R^2$. A continuous map $f : \R^2 \rightarrow \R^2$ is said to be a {\it topological horseshoe} for $N_1$ and $N_2$ if 
	\begin{eqnarray*}
		N_1 \stackrel{f} \Longrightarrow N_1 \, , \quad N_1 \stackrel{f} \Longrightarrow N_2 \, , \quad
		N_2 \stackrel{f} \Longrightarrow N_1 \, , \quad N_2 \stackrel{f} \Longrightarrow N_2 \, .
	\end{eqnarray*}
\end{definition}
\noindent Topological horseshoes are associated to symbolic dynamics as presented in 
Theorem 2 in~\cite{aGi19} and Theorem 18 in~\cite{ZGLICZYNSKI200432}, where the authors show that the existence of a horseshoe for a map $f$ provides a semi--conjugacy between $f$ and a shift map $\{0,1\}^{\Z}$, meaning that for any sequence of symbols $0$ and $1$ there exists an orbit generated by $f$ passing through the sets $N_1$ and $N_2$ in the order given by the sequence, guaranteeing the existence of ``any kind of orbit'' (periodic orbits, chaotic orbits, etc.).   \\
From the Definition~\ref{def_covering}, the covering relation $N_{1} \stackrel{P} \Longrightarrow N_{2}$ is verified if the three following conditions are satisfied:
\begin{enumerate}
	\item  the image $P(N_1)$ of $N_1$ lies in the strip between the top and the bottom edges of $N_2$,
	\item  the image of the left part of $N_1^-$ lies on the left of $N_2$,
	\item  the image of the right part of $N_1^-$ lies on the right of $N_2$;
\end{enumerate}
the conditions can be easily checked in Figure~\ref{fig:horseshoe} and, then, in Figure~\ref{fig:Q1Q2horseshoe}.

\subsection{Existence of a topological horseshoe}  
\begin{figure}
\includegraphics[height=9cm, width=10cm]{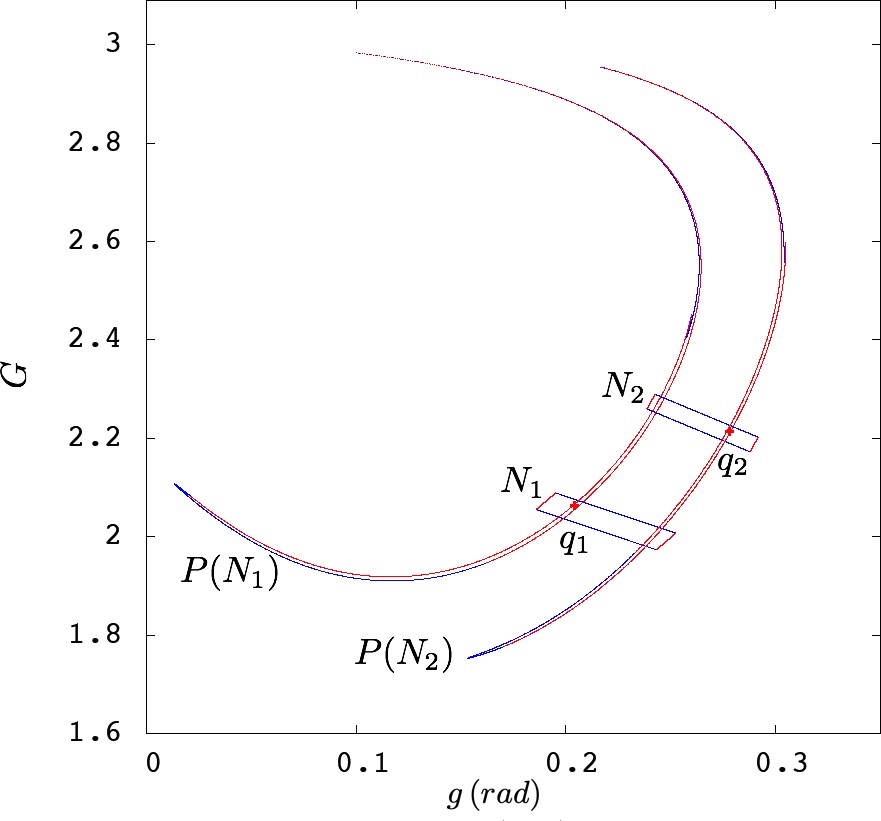}
\caption{\label{fig:Q1Q2horseshoe}
Horseshoe connecting the points $q_1$ and $q_2$ proving symbolic dynamics for the map $P$. Red represents the entry sets and their images and blue the exit sets and their images.
}
\end{figure}

In this section we describe how we construct explicitly a topological horseshoe for the Poincar\'e map of the Hamiltonian~\eqref{eq:AvH}. \\
We start by considering  one hyperbolic fixed point $q$ for the Poincar\'e map $P$ and we denote by $v^{s}$ and $v^{u}$, respectively, the stable and the unstable eigenvectors related to $DP(q)$. We construct a parallelogram $N$ containing $q$ whose edges are parallel to $v^{s}$ and $v^{u}$ and thus we define $N$ as 
\begin{eqnarray*}
	N = q + A v^{s} +B v^{u},
\end{eqnarray*} 
where $A$ and $B$ are suitable chosen closed real intervals. If the intervals $A$ and $B$ are sufficiently small, under the action of the map $P$, the parallelogram $N$ will be contracted in the stable direction and expanded in the unstable direction. We denote by $P(N)$ the image of $N$ through the map $P$. In practice, we choose two hyperbolic fixed points $q_{1}$ and $q_{2}$ having the important property of transversal intersection of their stable and unstable manifolds as shown in Figure~\ref{figJ:FPMan}. This property is a good indication of the existence of a topological horseshoe. Based on this couple of fixed points whose coordinates read 
\begin{eqnarray}
\left\{
\begin{aligned}
&q_{1}=(g_{1},G_{1}) = (0.203945459,2.06302430), \\ \notag 
&q_{2}=(g_{2},G_{2}) = (0.278077917,2.21418596),
\end{aligned}
\right.
\end{eqnarray}
we define  two sets $N_{1}, N_{2} \subset \R^2$ which are supports of two $h$--sets as follows:
\begin{eqnarray}
\left\{
\begin{aligned}
&N_{1} = q_1 + A_1 v_1^s +B_1 v_1^u, \\ \notag 
&N_{2} = q_2 + A_2 v_2^s +B_2 v_2^u,
\end{aligned}
\right.
\end{eqnarray}
where 
\begin{eqnarray}
\left\{
\begin{aligned}
&A_{1} = [-0.02,0.08] 	\subset \R, \hspace{0.35cm} \quad  B_{1} = [-0.025,0.01] \subset \R,  \\ \notag
&A_{2} = [-0.075,0.025]	 \subset \R, \quad B_{2} = [-0.02,0.01] \subset \R,
\end{aligned}
\right.
\end{eqnarray}
and $v_{1}^s,v_1^u,v_2^s,v_2^u$ are the stable and the unstable eigenvectors related to $q_{1},q_{2}$, respectively. Then the following covering relations hold
\begin{eqnarray*}
	N_1 \stackrel{P} \Longrightarrow N_1 \, , \quad N_1 \stackrel{P}\Longrightarrow N_2 \, , \quad
	N_2 \stackrel{P} \Longrightarrow N_1 \, , \quad N_2 \stackrel{P} \Longrightarrow N_2  \, ,
\end{eqnarray*}
proving the existence of a topological horseshoe for $P$, \ie existence of symbolic dynamics for $P$. The obtained horseshoe associated to $q_{1}$ and $q_{2}$ with the aforementioned parameters is illustrated in Figure\,\ref{fig:Q1Q2horseshoe}, providing the existence of symbolic dynamics. \\

\section{Conclusions and open problems}
\noindent This work originates from~\cite{gPi19}, where  it has been pointed out that the average $\UU_+$~\eqref{averaged} of the Newtonian potential  with respect to one of the two mean anomalies is an integrable function which in turn may be written as a function  of another function $t_+$, whose dynamics is completely known. The functional dependence (\ref{identity}) between these two functions, holding in the case of the planar problem,  has been pointed out in~\cite[Section 3]{gPi20b}. The identity (\ref{identity}) raises the very natural question whether and at which extent such relation has a consequence on the dynamics of the three--body problem. Giving an answer to such question is in fact  demanding, as 
it requires to  understand whether it is possible to find a region of phase space where the three--body Hamiltonian is well represented by its simple average (here ``simple average'' is used as opposite to ``double average'', most often encountered in the literature, e.g.,~\cite{laskarR95}) and, simultaneously, the kinetic term $\KK$ in~\eqref{KU} does not interfere with $\UU$ too much. In
\cite{gPi20b} it has been proved that if the total angular momentum $C$ of the system vanishes, by symmetry reasons, and using a well--suited perturbation theory, the librational motions of $t_+$ reported in Figure \ref{unperturbed phase portrait} (left) have a continuation in the averaged three--body problem. 
In this paper we  investigated the case $C\ne 0$. 
With purely pioneering spirit, in order to simplify the analysis, we focused on the very  peculiar situation where the two minor bodies have equal mass and we fixed an energy level once forever. We  believe that both such choices can be removed without affecting the results too much, because, as informally  discussed in the introduction, what really matters 
is the relative weight of $\KK$ and $\UU$. 
Figures \ref{unperturbed phase portrait} and \ref{fig:fig2J} not only show that, in our simplified model, this continuation is numerically evident, but also  exhibit the onset of chaos in certain zones, clearly highlighted along the paper using techniques of~\cite{aGi19}. 
Even though the results are encouraging, many questions are still pending (some of them have been  pointed out  in~\cite{gPi20b}), and we aim to face them in the future:
\begin{itemize}
\item[$({\rm Q}_1)$] If $C\neq 0$, is there a choice of parameters and phase space where the phase portrait of $\FF$ includes only librational motions?
\item[$({\rm Q}_2)$] In the case that the orbit of the planet is inner to the one of the asteroids, the phase portrait of $\UU_+$ includes a saddle and a separatrix through it (see~\cite[Figures 1, 2 and 3]{gPi20b}). How does this affect the three--body problem motions?
\item[$({\rm Q}_3)$] By~\cite{gPi19}, relation~\eqref{identity} has a generalisation to the spatial problem. What are the consequences on the spatial three--body problem? 
\item[$({\rm Q}_4)$] Is the onset of chaos in the averaged problem present also in the full (non--secular) system?
\item[$({\rm Q}_5)$] What can we prove analytically?
\item[$({\rm Q}_6)$] What can we prove with computer--assisted techniques?
\end{itemize}

\section*{Acknowledgments}
\noindent We are grateful to the anonymous reviewers for their stimulating remarks. We are indebted to C. Efthymiopoulos for a highlighting  discussion about how to control errors (Section~\ref{errors}) and to M. Guzzo for sharing his expertise on FLIs. We heartily thank U. Locatelli for an interesting talk during the meeting I-Celmech, that held in Milan, in February 2020. Figure~\ref{unperturbed phase portrait} has been produced using the software {\sc mathematica}\textsuperscript{\textregistered}.

\begin{appendices}

\section{The Hamiltonian}\label{Jacobi}
\noindent The impulses $p_0$, $y$, $y'$ conjugated to $r_0$, $x$, $x'$ in (\ref{x0}) are 
\begin{eqnarray}\label{impulsesp}p_0=y_0+y_1+y_2\ ,\quad y=y_1+\frac{\mu}{1+\mu}y_2-\frac{\mu}{1+\mu}p_0\ ,\quad y'=y_2-\frac{\kappa }{1+\mu+\kappa }p_0\ .
\end{eqnarray} 
If $r_0\equiv0\equiv p_0$, the transformation of coordinates defined by (\ref{x0}) and (\ref{impulsesp}) reduces to the injection 
\begin{eqnarray*}\left\{
\begin{array}{lll}
 x_0=-\frac{\mu}{1+\mu}x-\frac{\kappa }{1+\mu+\kappa }x'\\
x_1=\frac{1}{1+\mu}x-\frac{\kappa }{1+\mu+\kappa }x'\\
 x_2=\frac{1+\mu}{1+\mu+\kappa }x'
 \end{array}
 \right. \, , \qquad \left\{
 \begin{array}{lll}
 y_0=-y-\frac{1}{1+\mu}y'\\
y_1=y-\frac{\mu}{1+\mu}y'\\
 y_2=y'
 \end{array}
 \right.
 \end{eqnarray*}
 and the Hamiltonian~\eqref{first_Ham} becomes
\begin{eqnarray*} 
	\HH = 
	\frac{1 + \mu}{2 \mu m_0} {\vert y \vert}^2 +
	\frac{1 + \mu+\kappa }{2 (1+\mu) \kappa  m_0} {\vert y' \vert}^2 - 
	\frac{\mu m_{0}^2}{\vert x \vert}- \frac{\kappa  m_{0}^2}{\vert x' +
	\frac{\mu}{\mu+1} x\vert} 
	-
	\frac{\mu \kappa  m_{0}^2}{\vert x' -\frac{1}{\mu+1} x\vert}.
\end{eqnarray*}
Rescaling the coordinates via 
\begin{eqnarray*}
x\to (1+\mu) x\ , \quad y\to \frac{\mu}{1+\mu}  y\ , \quad  x'\to \beta^{-1}  x'\ , \quad y'\to \mu \beta y', 
\end{eqnarray*}
with $\beta$ as in (\ref{beta}) and multiplying the Hamiltonian $\HH$ by $(1+\mu)/\mu$, we obtain  $\HH$ as in (\ref{ham_rescaled}).
%

\section{Numerical setups and results}\label{App:NS}

\subsection{Choice of the parameters}\label{Choice of the parameters}
The analysis we have done is related to the choice of parameters and initial data we started with. The Hamiltonian ~\eqref{eq:AvH} is composed by three parts  
\begin{eqnarray*}
\HH =\HH_0+\sigma\KK+\sigma\UU=:\HH_0+\PP,
\end{eqnarray*}
 where the first one is the unperturbed and constant part depending on $\Lambda$, the second one represents the kinetic part and the third is the perturbing part. To ensure the non--resonant terms of $\PP$ to be small with respect to $\HH_0$  we choose, as mentioned in the introduction,
\begin{eqnarray}
\left\{
\begin{aligned}
&m_{0}= 1, \\ \notag
&\beta = 40, \\  \notag
&C= 75.597 \\ \notag 
&\Lambda = 3.099.  \notag
\end{aligned}
\right.
\end{eqnarray}
 The initial datum is taken to be
\begin{eqnarray}
\left\{
\begin{aligned}
& G_{\star} = -2.4915, \\ \notag
& R_{\star} = -0.0039, \\   \notag
& g_{\star}  =  1.4524,\\ \notag 
& r_{\star}  = 3132.069 . \notag
\end{aligned}
\right.
\end{eqnarray}
Note that $R_\star$ and $r_\star$ verify (\ref{circular}) but are not exactly centred at $0$ and $r_0$ because the $r$--component of the Hamiltonian vector--field vanishes for $R=0$, while it needs to be different from zero in order that the Poincar\'e map is well defined. The values of $G_\star$ and $g_\star$ have been empirically chosen such that the orbit from from $(G_\star, R_\star, g_\star, r_\star)$ is  approximately periodic  and hence the Poincar\'e map is well defined. 

\subsection{Flow}\label{errors}
The Hamiltonian equations of motion have been numerically propagated using a fixed time--step \textsc{RK4} method~\cite{pr92}. Even though the step has been kept fixed, no numerical issues have been encountered and the integration times were reasonable for the whole numeric exploration.\\
 Under the choice of our time--step $\delta$, the flow--map preserves the Hamiltonian itself, a conserved quantity (\textit{first integral}), with a relative error of about $10^{-14}$ for stable orbits and $10^{-12}$ for chaotic orbits on a arc length of about $\tau \sim 10^2$ orbital revolutions. Besides the first integral being numerically well preserved, the quality of the integration has been assessed further using a forwards/backwards strategy. The method consists in propagating forwards in time (say on $[0,\tau]$) the Cauchy problem 
\begin{eqnarray}
\left\{
	\begin{aligned}
		&\dot{x}=v_{\mathcal{H}}(x), \\ \notag
		&x(0)=x_{0},
	\end{aligned}
\right.
\end{eqnarray} 
 and then to back--propagate (from $\tau$ to 0) the new Cauchy problem 
\begin{eqnarray}
\left\{
	\begin{aligned}
		&\dot{x}=v_{\mathcal{H}}(x), \\ \notag
		&x(\tau)=x_{\tau}
		\end{aligned}
\right.
\end{eqnarray} 
 where the initial seed $x_{\tau}$ is obtained from the forward numerical flow--map, $x_{\tau}=\Phi^{\tau}(x_{0})$. Then the relative error
 \begin{eqnarray*}
 	\Delta=\frac{\norm{x_{0}-\Phi^{-\tau}\big(\Phi^{\tau}(x_{0})\big)}}{\norm{x_{0}}}
 \end{eqnarray*} 
 is estimated. On a selection of orbits, we found $\Delta$ to be of the order of $10^{-12}$ for regular orbits, $10^{-8}$ for chaotic orbits on timescale of about $10^{2}$ orbital revolutions. 

\subsection{Poincar\'e mapping $P$}
The construction of the Poincar\'e map $P$ is based on the time evolution of the whole flow and a bisection procedure. Given an initial point $z$, to find its next state $z'=P(z)$ we compute $x(t)=\Phi^{t}(x_{0})$, $x=l(z)$, until following conditions are met:
\begin{enumerate}
	\item \textit{Section condition}: $X=(x_{1}(t),x_{2}(t),x_{4}(t)) \in \Sigma$ up to a numerical tolerance $\varepsilon_{\Sigma}=10^{-10}$. This step relies on a bisection method halving the length of the numerical step $\delta$ until we drop under the tolerance $\varepsilon_{\Sigma}$. 
	\item \textit{Orientation condition}: The scalar product $\dot X(0)\cdot  \dot X(t)$ is positive, meaning that the orbit is intersecting the plan $\Sigma$ in the same direction as the starting point.
	\item \textit{First-return condition}: for $\tau<t$, neither (1) and (2) are fulfilled.
\end{enumerate}

\subsection{Coordinates of the fixed--points of $P$}\label{sub:CoordFP}
Below we provide the coordinates of the fixed--points of $P$ (periodic orbits of $\mathcal{H}$). 

\begin{center}
\begin{varwidth}{\linewidth}
\begin{verbatim}
# Coordinates of the elliptic fixed points                                                                                                                                                                                                    
#  	     G            g (rad)                                                                                                                                                                                                                       
#-------------------------------------------                                                                                                                                                                                                    
      -2.49155       1.45245
      -1.04685       1.73094
      -2.91949       1.95066
\end{verbatim}
\end{varwidth}
\end{center}
\medskip
\begin{center}
\begin{varwidth}{\linewidth}
\begin{verbatim}
# Coordinates of the hyperbolic fixed points                                                                                                                                                                                                    
#  	     G            g (rad)                                                                                                                                                                                                                        
#-------------------------------------------                                                                                                                                                                                                    
       2.06302       0.20395
       2.21419       0.27808
       0.03851       0.33259
       2.47589       0.34655
       2.81488       0.40502
       3.04924       0.43647
       3.09865       0.44249
      -2.84323       0.55513
       2.75151       0.57177
       3.05336       0.58816
       3.09883       0.59055
      -2.61168       1.35169
       2.76024       1.61321
       2.68138       2.39082
       2.52039       2.51911
       2.31386       2.60074
       2.49651       2.61696
       1.85433       2.62309
       1.75010       2.62341
       2.43689       2.75722
       2.33537       2.90395
       2.22839       3.01548
\end{verbatim}
\end{varwidth}
\end{center}

\medskip

\section{The Fast Lyapunov Indicator \& dynamical timescales}\label{App:FLI}
\begin{figure}
\includegraphics[height=5cm, width=13cm]{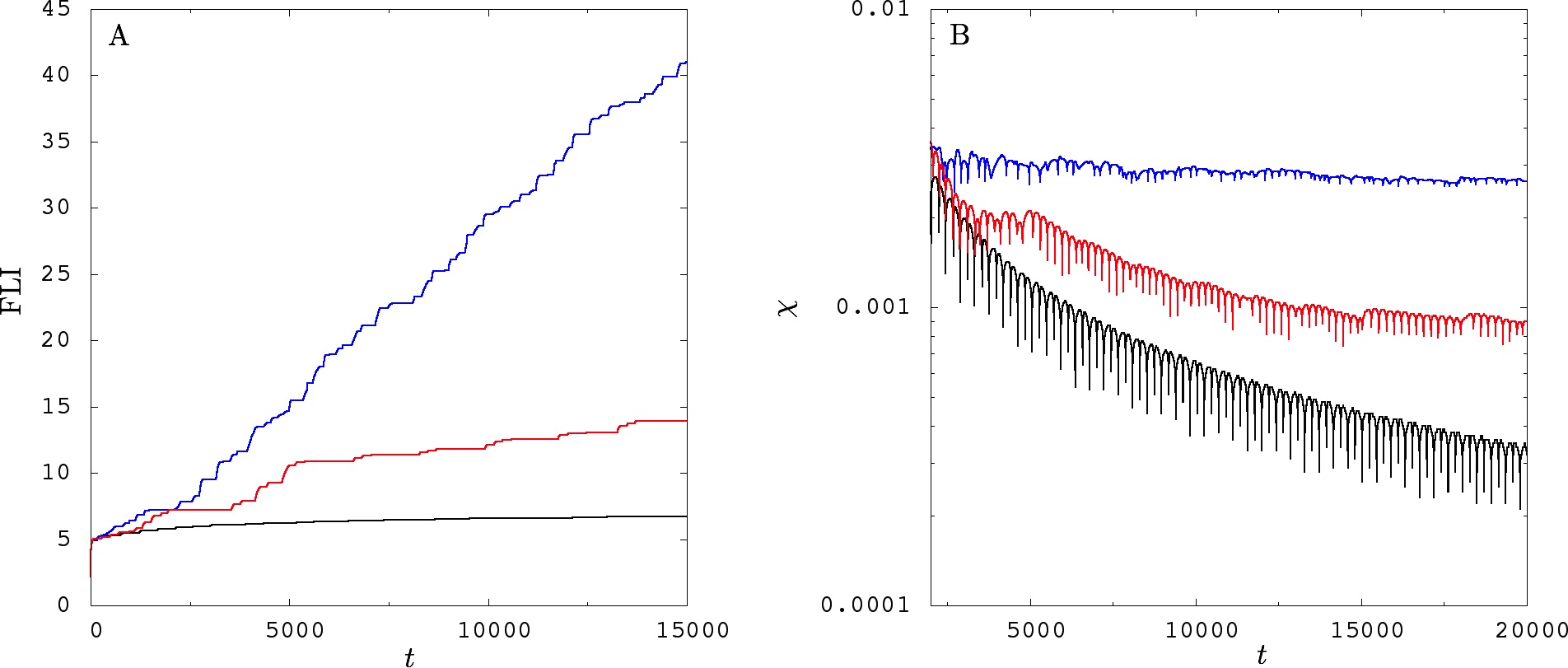}
\caption{\label{figJ:FLIcalibration}
On the left: calibration of the finite time chaos indicators (FLI) for three distinct orbits. After a transient time of about $t \sim 5,000$ (\ie $\sim 10$ orbital revolutions) the discrimination of the nature of the orbits is  sharp enough. On the right:  time evolution of the maximal Lyapunov exponents $\chi$. For chaotic orbits, they define Lyapunov timescales of about $0.76$ orbital revolutions.  
}
\end{figure}

\noindent The Fast Lyapunov Indicator (FLI) is an easily implementable tool suited to detect phase--space structures and local divergence of nearby orbits.  It has a long--lasting tradition with problems motivated by Celestial Mechanics~\cite{cFr97}. The indicator can be used in the context of deterministic  ODEs, mappings, and is able to detect manifolds and global phase--space structures~\cite{cFr00,mGu14,eLe16}. A large literature exists  with the FLI tested on idealised systems (\eg low dimensional quasi--integrable Hamiltonian system~\cite{cFr00,mGu13}, drift in volume--preserving mappings~\cite{nGu17}) but  also on many applied gravitational problems across a variety of scales, ranging from the near--Earth space environment~\cite{jDa18} to exoplanetary systems~\cite{riPa15}. For simplicity, let us  present the tool in the case of ODEs. Let us assume we are dealing with a $n$--dimensional autonomous ODE system. If the system is non--autonomous, we classically extend the dimension of the phase--space by $1$ dimension. The FLI indicator is based on the \textit{variational dynamics} in $\mathbb{R}^{2n}$,
\begin{eqnarray}
\left\{
\begin{aligned}
	&\dot{x}=f(x),\notag\\
	&\dot{w}=Df(x)\cdot w,
\end{aligned}
\right.   
\end{eqnarray}
$w \in T_{x}M$,
and is defined at time $t$ as 
\begin{eqnarray}\label{FLI}
	\textrm{FLI}(x_{0},w_{0},t) = {\rm sup}_{0 \le \tau \le t} \log \norm{w(\tau)}.
\end{eqnarray}
The FLI is able to distinguish quickly the nature of the orbit emanating from $x_{0}$. Orbits containing the germ of hyperbolicity will have their final FLI values larger than regular orbit (for the same horizon time $\tau$). More precisely, chaotic orbits will display a linear growth (with respect to time) of their FLIs, whilst regular orbits have their FLIs growing logarithmically. In order to reduce the parametric dependence of the FLIs upon the choice of the initial tangent vector, the FLIs are computed over an orthonormal basis of the tangent space (\ie we compute Eq.~\eqref{FLI} 4 times with a different initial $w_0$) and averaged~\cite{mGu13}. As a rule of thumb, the FLI is computed over a few Lyapunov times $\tau_{\mathcal{L}}$, but in practice, the choice of the final $\tau$ requires a calibration procedure by  testing few orbits.  By computing FLIs on discretised domains of initial conditions, the color coding of the FLIs (using a divergent color palette) reveals the global topology of the phase--space (\eg web of resonances and preferred routes of transport, see~\cite{mGu13}) furnishing a so--called \textit{stability map}. The Lyapunov time  $\tau_{\mathcal{L}}$ is obtained as the inverse of the \textit{maximal Lyapunov characteristic exponent} (we refer to~\cite{chSk10} for computational aspects related to characteristic exponents),
\begin{eqnarray*}
	\tau_{\mathcal{L}} = 1/\chi,
\end{eqnarray*}
 where $\chi$ denotes the maximal Lyapunov characteristic exponent
 \begin{eqnarray*}
 	\chi(x_{0},w_{0}) = \lim_{t \to +\infty}\frac{1}{t} \log \norm{w(t)}.
 \end{eqnarray*}
Stable orbits do satisfy $\chi \to 0$ and hence $\tau_{\mathcal{L}}$ tends to be large. On the contrary, chaotic orbits are characterised  by $\chi \to r \in \mathbb{R}^{\star}_{+}$ and therefore  $\tau_{\mathcal{L}}$ converges to a finite value. The panel shown in Figure\,\ref{figJ:FLIcalibration} presents the calibration procedure based on three  orbits. The stable orbit displayed in black (logarithmic growth of the FLI) admits for initial condition $(G,g)=(-2,\pi)$. The two others orbits are chaotic but one (red) is less hyperbolic than the other (blue). The respective initial conditions read $(G,g)=(-2,1.6)$ and $(G,g)=(2,\pi)$. As it is observed, after a transient time of about $t \sim 5,000$ (\ie $10$ orbital revolutions), safe conclusions can be formulated regarding the stability of the orbits (left panel). The respective maximal Lyapunov characteristic exponents are presented in the right panel of Figure\,\ref{figJ:FLIcalibration}.  The inverse, the Lyapunov time, defines timescales of $\sim 380$ for the most chaotic one (which is about $0.76$ revolutions) and $\sim 1,100$ for the second chaotic one ($2.2$ orbital revolutions).   

\end{appendices}

\bibliographystyle{apalike}
\bibliography{biblio}

\end{document}